\documentclass{elsart}
\usepackage{etex}

\usepackage[numbers]{natbib}
\usepackage{easymat}
\usepackage{easybmat}

\usepackage{mathtools}
\usepackage{pslatex,times}
\usepackage{amssymb}
\usepackage{mathpazo}
\usepackage{easyeqn}
\usepackage[leqno,fleqn,intlimits]{empheq}

\usepackage[T1]{fontenc}
\usepackage{graphicx}
\usepackage{psfrag}
\usepackage{wrapfig}
\usepackage{multirow,tabls}
\usepackage{booktabs}
\usepackage{psfrag}

\usepackage[a4paper,textwidth=16cm, textheight=23cm,
					  hmarginratio=1:1,vmarginratio=1:1]{geometry}

\newcommand{\dt}{\delta t}
\newcommand{\pd}{\partial}

\newcommand{\bu}{{\mathbf u}}

\newcommand{\bff}{{\mathbf f}}
\newcommand{\bF}{{\mathbf F}}
\newcommand{\bg}{{\mathbf g}}

\newcommand{\bB}{{\mathbf B}}

\newcommand{\bA}{{\mathbf A}}

\newcommand{\bx}{{\mathbf x}}

\newcommand{\lap}{\Delta}
\newcommand{\source}{S}
\newcommand{\bsource}{{\mathbf S}}

\newcommand{\bnabla}{\boldsymbol{\nabla}} 
\newcommand{\blap}{\boldsymbol{\Delta}}
\newcommand{\grad}{{\nabla}}
\newcommand{\bgrad}{{\bnabla}}
\newcommand{\curl}{{\nabla} \times}
\newcommand{\divv}{{\nabla} \cdot}

\newcommand{\bn}{\hat{\mathbf n}}
\newcommand{\be}{\hat{\mathbf e}}
\newcommand{\er}{\hat{\mathbf e}_r}

\newcommand{\ez}{\hat{\mathbf e}_z}
\newcommand{\et}{\hat{\mathbf e}_\theta}

\newcommand{\mx}[1]{\left[{#1}\right]}

\newcommand{\halfM}{\left\lfloor\frac{M}{2}\right\rfloor}

\newcommand{\pdt}{\partial_t}
\newcommand{\pdr}{\partial_r}
\newcommand{\pdrz}{\partial_{rz}}
\newcommand{\pdth}{\partial_\theta}
\newcommand{\pdthz}{\partial_{\theta z}}
\newcommand{\pdz}{\partial_z}

\newcommand{\pmh}{\pm\frac{h}{2}}

\newcommand{\mF}{\mathcal{F}}

\newcommand{\mT}{\mathcal{T}}

\newcommand{\mQ}{\mathcal{Q}}
\newcommand{\ir}{\frac{1}{r}}
\newcommand{\Rm}{Rm}
\newcommand{\iR}{Re^{-1}}
\newcommand{\iRm}{\Rm^{-1}}
\newcommand{\laph}{\lap_h}
\newcommand{\gradh}{\grad_h}
\newcommand{\divh}{\grad_h\cdot}
\newcommand{\ie}{i.e.}
\newcommand{\eg}{e.g.}
\newcommand{\p}{^\prime}
\newcommand{\domain}{\Omega}
\newcommand{\angvel}{\omega}
\newcommand{\bdy}{{\pd\Omega}}
\newcommand{\infmat}{\mathcal{C}}
\newcommand{\truecond}{c}

 \numberwithin{equation}{section}
\begin{document}
\begin{frontmatter}
\title{Poloidal-toroidal decomposition in a finite cylinder.
I. Influence matrices for the magnetohydrodynamic equations}
\author{Piotr Boronski and Laurette S. Tuckerman}
\address{LIMSI-CNRS, BP 133, 91403 Orsay, France}
\begin{abstract}
 The Navier-Stokes equations and magnetohydrodynamics equations are written in terms of poloidal and toroidal potentials in a finite cylinder. This formulation insures that the velocity and magnetic fields are divergence-free by construction, but leads to systems of partial differential equations of higher order, whose boundary conditions are coupled. The influence matrix technique is used to transform these systems into decoupled parabolic and elliptic problems. 
The magnetic field in the induction equation is matched to that in an exterior vacuum by means of the Dirichlet-to-Neumann mapping, thus eliminating the need to discretize the exterior.
The influence matrix is scaled in order to attain an acceptable condition number.
\end{abstract}
\end{frontmatter}

\section{Motivation and Governing Equations}
The requirement that velocity and magnetic fields be solenoidal,
i.e.~divergence-free, represents one of the most challenging difficulties
in hydrodynamics and in magnetohydrodynamics 
\cite{Marques90,Marques93,Tuckerman89,Rempfer06,LopezMarquesShen02,Brackbill,Chan01}.
For the velocity field, this condition is the fundamental approximation used
in incompressible fluid dynamics.  For the magnetic field, this condition is
the statement of the non-existence of magnetic monopoles.

Two main approaches exist for imposing this requirement.  The first 
is to use three field components and to project three-dimensional fields
onto a divergence-free field.  In an incompressible fluid, the pressure serves
to counterbalance the nonlinear term which is the source of the divergence in
the Navier-Stokes equations; the pressure also plays this role numerically.
The divergence of the Navier-Stokes equations is taken, leading to a Poisson
problem for the pressure. However, the boundary conditions on the equations
for $(\bu,p)$ involve only the velocity, leading to coupling between the
equations to be solved for $\bu$ and $p$ \cite{Tuckerman89,Rempfer06}.  The
coupled equations can be solved in several stages by a Green's function or
influence matrix method \cite{Tuckerman89}.  In projection-diffusion schemes,
approximate boundary conditions are imposed for the pressure
\cite{LopezMarquesShen02}.  For magnetic fields, however, the exact
evolution of the equations conserves divergence and there exists no analogue
to the pressure.  Thus if the numerical algorithm creates divergence, there is
no mechanism for eliminating it and it may accumulate \cite{Brackbill}.  For
this reason, magnetohydrodynamic codes sometimes include a fictitious term
analogous to the hydrodynamic pressure, which must be treated numerically
\cite{Chan01}.

The second approach, which is the focus of this paper, 
is to express fields in such a way that they
are divergence-free by construction.
It can be proved that a field $\bF$ which is solenoidal (divergence-free)
in a simply connected domain can be written as:
\begin{equation}
\bF=\curl\left(\psi \be\right)+\curl\curl\left(\phi\be\right)
\label{eq:poltor}\end{equation}
where $\be$ denotes a unit vector. 
In addition to being divergence-free, $\bF$ has the 
advantage of involving only two scalar fields. 
This makes more economical use of computer memory and 
allows all calculations to be implemented using only scalar fields.

Equations governing the evolution of the two potentials are derived by taking
the curl and double curl of the original equations, increasing the order of
the differential equations.  In addition, boundary conditions, some also of
high order, couple the two potentials.  In certain geometries with two
periodic directions, 
these are only minor obstacles \cite{Marques90}.  In spectral treatments of
such geometries, the basis functions insure periodicity, which is preserved
under differentiation and addition. At most, special consideration must be
given to constant modes.  The standard examples are a spherical geometry
\cite{Marcus81,Glatzmaier84,Dudley89,Glatzmaier95,Tilgner97,Hollerbach00} or a
three-dimensional Cartesian geometry with one bounded direction and two
perpendicular periodic directions, such as channel flow \cite{Squire,Schmid}.
Other applications are in a cylindrical geometry with periodic $z$ and $\theta$
directions \cite{Marques90,Antonijoan98,Willis02}.

In geometries with more than one nonperiodic direction, far more care is
required.  Marques \cite{Marques90} gave a detailed analysis of the
poloidal-toroidal decomposition for the Navier-Stokes equations and its
formulation and validity for general topologies.  This analysis was then put
into practice in a linear stability analysis of Rayleigh-B\'enard convection
in a finite cylindrical geometry \cite{Marques93}. However, the governing
equations derived in \cite{Marques93} contain large linear systems that couple
the potentials and their laplacians and bilaplacians, but whose solution would
be required in implicit time integration.  Analogous problems arise
in the other formulations of incompressible fluid dynamics.
In the 2D streamfunction-vorticity formulation, the equations
for the vorticity and the streamfunction are coupled by the
fact that boundary conditions exist only for the streamfunction and none
on the vorticity. In the $(\bu,p)$ primitive variable formulation, 
the pressure is the solution to a Poisson problem for
which the appropriate boundary condition is that 
the velocity be divergence-free \cite{Tuckerman89,Rempfer06}.

Our primary goal in this paper is to demonstrate that the high-order equations
can be separated via the influence matrix technique into a sequence of
problems of lower order, each with its own boundary conditions, as was done
for the primitive variable formulation in \cite{Tuckerman89}.  This makes
implicit time integration feasible for the poloidal-toroidal decomposition in
geometries with two non-periodic directions. A secondary goal is to carry out
the same analysis for a magnetic field which is governed by the induction
equation and which generalizes the Navier-Stokes equation by the inclusion of
the Lorentz force.

The equations we will consider are the magnetohydrodynamic
equations:
\vspace*{-0.5cm}
\begin{subequations}
\label{eq:varr}
\begin{align}
\pdt \bu+(\bu\cdot \bnabla) \bu &= (\bB \cdot \bnabla) \bB + \iR \lap \bu - \bnabla (p+\frac{B^2}{2}) \label{v1}\\
\bnabla \cdot \bu &= 0 \label{v2}
\end{align}
\end{subequations}
\vspace*{-1cm}
\begin{subequations}
\label{eq:barr}
\begin{align}
\pdt \bB&=\curl (\bu\times \bB)+ \iRm\lap \bB \label{b1}\\
\nabla \cdot \bB &= 0 \label{b2}
\end{align}
\end{subequations}
where $Re$ is the usual hydrodynamic Reynolds number and 
$\Rm$ the magnetic Reynolds number.
Equations (\ref{eq:varr}) and (\ref{eq:barr}) are of different types:
for a divergence-free magnetic field $\bB$, all the terms of 
(\ref{eq:barr}) have zero divergence as well, but this is not the case 
for (\ref{eq:varr}).

The velocity and magnetic fields are to be calculated in a finite cylinder.
We consider specifically the case in which the flow is driven
by rotating upper and lower disks, although our method does not
depend on this. For disks rotating in opposite directions
this configuration is called the von K\'arm\'an flow
\cite{vonKarman,VKS,Nore03}.  The magnetic field inside the cylinder is
required to match the field outside, which goes to zero at infinity.  These
boundary conditions are expressed as:
\vspace*{-0.5cm}
\begin{subequations}
\label{eq:cond-u}
\begin{align}
\bu&= 0 && \text{at $r=1$} , \label{eq:cond-uB1}\\
\bu&= r\angvel_\pm\et && \text{at $z=\pmh$},\label{eq:cond-uB2}
\end{align}
\end{subequations}
\vspace*{-1cm}
\begin{subequations}
\label{eq:cond-B}
\begin{align}
\bB^{int}-\bB^{ext} &= 0 && \text{on $\bdy$},\\
\bB&=0 && \mbox{at infinity}.\label{eq:cond-uB4}
\end{align}
\end{subequations}
where $\Omega$ denotes the interior domain (the cylinder) and $\bdy$ is its boundary. 

\begin{figure}[h]
\centering
\includegraphics[width=4cm]{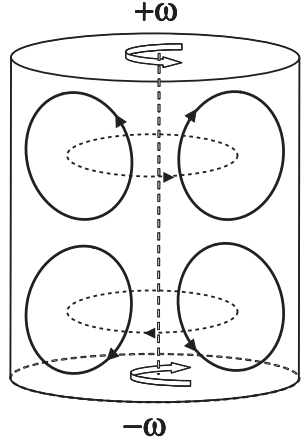}
\caption{Axisymmetric flow between counter-rotating disks.
Poloidal component: solid curves.
Toroidal component: dashed curves.}
\label{fig:polo-toro}
\end{figure}

The poloidal and toroidal components for this configuration in the
axisymmetric case with $\be=\ez$ are illustrated in figure
\ref{fig:polo-toro}.  The toroidal flow corresponds to motion with only
azimuthal velocity.  The poloidal flow forms recirculation rolls in the
$(r,z)$ plane.  For a non-axisymmetric flow, there is no clear correspondence
between each potential and a simple topological structure.

In section \ref{sec:decomposition}, we give a general description of the
poloidal-toroidal decomposition.  In section \ref{sec:hydro_BC}, we then
specialize to the Navier-Stokes equations in a finite cylinder, formulating
the boundary conditions for this case.  In section \ref{sec:nested_u} we show
how to decouple the equations and boundary conditions via the influence matrix
technique.  Finally, in section \ref{sec:MagneticSolver}, we present the
equations and boundary conditions for the induction equation which governs the
magnetic field, and the corresponding influence matrix.

\section{Poloidal-toroidal decomposition}
\label{sec:decomposition}

\subsection{Governing equations}

\begin{figure}[t]
\centering
\includegraphics[height=6cm]{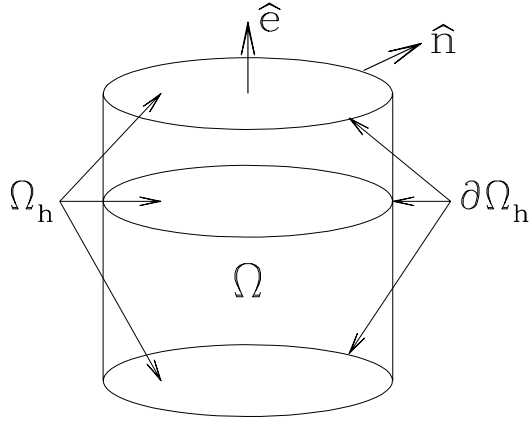}
\caption{Geometry for potential variable formulation. $\domain$ is a
cylindrical domain.  The vector $\be$ points in the distinguished vertical
direction, here $\be_z$.  $\domain_h$ are slices of $\domain$ perpendicular to
$\be$, here disks.  The boundary of $\domain_h$ is $\pd\domain_h$, here a
circle.  The vector $\bn$ is normal to both $\be$ and to $\domain_h$; here
$\bn = \be_r$.}
\label{fig:compat}
\end{figure}

The poloidal-toroidal decomposition generalizes to three
dimensions the two-di\-men\-sional stream\-function-vorticity formulation.
We follow the analysis and notation of \cite{Marques90}, 
but specializing to the case of a domain which is contractible 
to a point (\ie~has no holes). Then:
\begin{equation}
\divv\bF=0 \quad \Leftrightarrow \quad \bF=\curl {\mathbf A}
\label{eq:divcurl}\end{equation}
A distinguished direction and associated unit vector $\be$ is selected
and $\bA$ can be decomposed such that: 
\begin{equation}
\bF = \curl \psi \be+ \curl \curl \phi \be
\label{eq:poltor2}\end{equation}
The direction $\be$ is called vertical and those
perpendicular to $\be$ are called horizontal; 
see figure \ref{fig:compat}.
%
A number of possibilities exist for $\be$. Among these,
the choices $\be=\ez$ (in Cartesian or cylindrical coordinates)
or $\be=\be_\rho$ (the spherical radius)
decouple $\psi$ and $\phi$ in the diffusive operators since:
\begin{subequations}
\label{eq:curls}
\begin{alignat}{2}
\be\cdot\bF &= -\laph\phi, \qquad
&\be\cdot\lap\bF &= -\lap\laph\phi, \\
\be\cdot\curl\bF &= -\laph\psi, \qquad
&\be\cdot\curl\lap\bF &= -\lap\laph\psi, \\
\be\cdot\curl\curl\bF &= \lap\laph\phi,\qquad
&\be\cdot\curl\curl\lap\bF &= \lap\lap\laph\phi.
\end{alignat}
\end{subequations}
where $\laph$ is the two-dimensional Laplacian acting in the horizontal
directions, \ie, those perpendicular to $\be$.  (The decoupling
\eqref{eq:curls} does not hold \cite{Marques90} when the cylindrical radius
$\er$ is chosen as the distinguished direction $\be$).

The equations for the velocity potentials are derived by taking
the $\be$ component of the single and double curl of \eqref{eq:varr};
those for the magnetic potentials are derived by taking the
$\be$ component itself and the single curl of \eqref{eq:barr}.
The difference arises from the fact that all the terms of
(\ref{eq:barr}) are divergence-free and there is no pressure to eliminate.
Combining \eqref{eq:varr}-\eqref{eq:barr} and \eqref{eq:curls}
leads to the evolution equations for the scalar potentials:
\vspace*{-0.5cm}
\begin{subequations}
\label{eq:potMHD_u}
\begin{align}
(\pdt-\iR\lap)\laph\psi_u &= \be\cdot\curl \bsource_u
\label{eq:potMHD_u1}\\
(\pdt-\iR\lap)\lap\laph\phi_u &= -\be\cdot\curl\curl \bsource_u
\label{eq:potMHD_u2}
\end{align}
\end{subequations}
\vspace*{-1cm}
\begin{subequations}
\label{eq:potMHD_B}
\begin{align}
(\pdt-Rm^{-1}\lap)\laph\phi_B &= \be\cdot \bsource_B
\label{eq:potMHD_B1}\\
(\pdt-Rm^{-1}\lap)\laph \psi_B &= \be\cdot\curl \bsource_B
\label{eq:potMHD_B2}
\end{align}
\end{subequations} 
\vspace*{-0.5cm}
where:
\begin{subequations}
\begin{align}
\bsource_u &\equiv (\bu\cdot\bnabla)\bu-(\bB\cdot\bnabla)\bB\label{eq:def_su}\\
\bsource_B &\equiv -\curl(\bu\times\bB)\label{eq:def_sB}
\end{align}
\label{eq:def_su_sB}
\end{subequations}
Equations \eqref{eq:potMHD_u}-\eqref{eq:potMHD_B} are not all of the same
order in the vertical and horizontal directions.  For example, for the
velocity, \eqref{eq:potMHD_u1} is $2^{\rm nd}$ order in the vertical direction
and $4^{\rm th}$ order in the horizontal directions, while
\eqref{eq:potMHD_u2} is $4^{\rm th}$ order in the vertical direction and
$6^{\rm th}$ order in the horizontal directions.  A corresponding number of
boundary conditions are required for the velocity potentials, a total of
(2+4)/2=3 conditions at each vertical boundary and (4+6)/2=5 at each
horizontal boundary for $\bu$.  The conditions at the vertical boundaries are
those corresponding to the physical problem. At the horizontal boundaries, the
physical conditions must be supplemented by two additional conditions whose
derivation is the subject of the remainder of this section.

\subsection{Gauge freedom}
\label{sec:gauge_freedom}

The poloidal-toroidal formulation \eqref{eq:poltor2} contains a gauge freedom
for the choice of $\psi$ and $\phi$, which is identified by finding the class
of potentials satisfying the homogeneous problem $\bF=0$. For $\be=\ez$
(Cartesian or cylindrical coordinate) or $\be=\be_\rho$ (spherical radius),
this leads to:
\begin{subequations}
\begin{EQA}[rrll]
\bF^{hom} = 0=&\curl\left(\psi^{hom}\be\right)&+\curl\curl\left(\phi^{hom}\be\right)&\\
0=&\be\times\grad_h\psi^{hom}&+\grad_h\pd_e\phi^{hom}-\left(\lap_h\phi^{hom}\right)\be\\
&&\Downarrow&\\ 
\be\cdot\bF^{hom} =0  \quad \Rightarrow \quad &\lap_h\phi^{hom} &= 0&\label{eq:gauge_cond_a}\\
\be\times\bF^{hom}=0 \quad \Rightarrow \quad  &\be\times\grad_h\psi^{hom} &= -\grad_h\pd_e\phi^{hom}&\label{eq:gauge_cond_b}
\end{EQA}
\label{eq:gauge_cond}
\end{subequations}
where $e$ denotes the coordinate corresponding to $\be$.
The existence of $\psi^{hom}$ satisfying \eqref{eq:gauge_cond_b} for all
$\phi^{hom}$ satisfying \eqref{eq:gauge_cond_a} 
is demonstrated as follows. Condition 
\eqref{eq:gauge_cond_a} implies:
\begin{equation}
\grad_h\cdot(\grad_h\pd_e\phi^{hom})=0
\label{eq:gauge_cond_b_dz}
\end{equation}
Applying \eqref{eq:divcurl} to the simply-connected two-dimensional
domain slices perpendicular to $\be$, 
\eqref{eq:gauge_cond_b_dz} implies that there exists a $\psi^{hom}$ satisfying:
\begin{equation}
\grad_h\pd_e\phi^{hom}=\gradh\times(-\psi^{hom}\be)=\be\times\grad_h\psi^{hom}
\end{equation}
Thus, the poloidal potential $\phi$ is determined up
to a harmonic function on each domain slice perpendicular to $\be$:
\begin{subequations}
\begin{equation}
	\phi \sim \phi + \phi^{hom} \quad;\quad\laph \phi^{hom}=0
	\label{eq:phi_gauge}
\end{equation}
while $\psi$ is determined up to an arbitrary function of the
coordinate $e$:
\begin{equation}
	\psi \sim \psi + h(e)
	\label{eq:psi_gauge}
\end{equation}
\end{subequations}
The choice of gauge constitutes one of the two additional conditions required.

\subsection{Compatibility condition}
\label{sec:compat-cond}

We have not yet demonstrated the equivalence between the potential and
primitive variable formulations.  Since the curl of equations \eqref{eq:varr}
and \eqref{eq:barr} were taken, they gained an additional degree of freedom
which must be fixed in such a way that these equations in potential form
\eqref{eq:potMHD_u}-\eqref{eq:potMHD_B} define the same velocity $\bu$ and
magnetic field $\bB$ as the original MHD equations
\eqref{eq:varr}-\eqref{eq:barr}.  We will require the fact that on a
simply-connected domain, a field is a gradient if and only if it is curl-free:
\begin{equation}
{\mathbf f}=\grad p \quad \Leftrightarrow \quad \curl{\mathbf f}=0
\label{eq:curlgrad}
\end{equation}
which is a consequence of Stokes' theorem.
We will first write \eqref{eq:varr}-\eqref{eq:barr} 
in a compact form, which will let us use a
common form for \eqref{eq:potMHD_u} and \eqref{eq:potMHD_B}:
\vspace*{-0.5cm}
\begin{subequations}
\label{eq:compact_form}
\begin{eqnarray}
\bff_u\equiv & \left(\pdt-\iR\lap\right)\bu + \bsource_u 
& =-\grad\left(p+B^2/2\right)
\label{eq:def_fu}\\
\bg_u\equiv & \curl{\mathbf f}_u & = 0\label{eq:def_gu}\\
\bg_B\equiv & \left(\pdt-Rm^{-1}\lap\right)\bB + \bsource_B 
& = 0\label{eq:def_gB}
\end{eqnarray}
\end{subequations}%
where \eqref{eq:def_fu} and \eqref{eq:def_gu} are equivalent, by
\eqref{eq:curlgrad}.  Then we can write the primitive variable formulation
\eqref{eq:varr}-\eqref{eq:barr} and potential formulation
\eqref{eq:potMHD_u}-\eqref{eq:potMHD_B} using for either $\bg=\bg_u$ or
$\bg=\bg_B$:
\begin{equation}
\begin{array}{rcl}
\text{\it primitive variables}&&\text{\it potential formulation}\\
\bg=0&\Rightarrow& 
	\left\{
	\begin{array}{rl}
	\be\cdot\bg &= 0\\
	\be\cdot\curl\bg &= 0
	\end{array}
	\right.
\end{array}\label{eq:rightonly}
\end{equation}
Marques \cite{Marques90} proves that in a simply connected domain $\domain$, 
the potential and primitive 
variable formulations are equivalent if additional conditions are satisfied:
\begin{subequations}
\label{eq:equiv_systems}
\begin{empheq}[left={\bg=0\quad \Leftrightarrow \quad\empheqlbrace\quad}]{align}
\be\cdot\bg = 0& \quad \text{in $\domain$}\label{eq:eg0}\\
\be\cdot\curl\bg = 0& \quad \text{in $\domain$}\label{eq:ecg0}\\
\divv\bg =0& \quad \text{in $\domain$}\label{eq:comp_nodiv}\\
\bn\cdot\bg = 0& \quad \text{on $\pd\domain_h$}\label{eq:comp_norm_proj}
\end{empheq}
\end{subequations}
In \eqref{eq:comp_norm_proj}, $\bn$ is the vector normal to
the boundary $\pd\domain_h$ of slices perpendicular to $\be$. 
We recall that in our case, $\be=\ez$,
the slices $\domain_h$ are disks, their boundaries
$\pd\domain_h$ are circles, and $\bn=\er$ is the radial unit vector.
We illustrate this geometry in figure \ref{fig:compat}.

The rightwards implication of \eqref{eq:equiv_systems} is obvious. The
leftwards implication of \eqref{eq:equiv_systems} is proved as follows.
We first use \eqref{eq:eg0} and the two-dimensional Stokes' Theorem
\eqref{eq:curlgrad} to introduce a scalar function $\kappa$
\begin{equation}
\left.\begin{array}{r}
0 = \be\cdot \bg \\
0 = \be\cdot \curl \bg \end{array}\right\}\Rightarrow
\bg=\gradh \kappa
\end{equation}
recalling that the subscript $h$ restricts differential operators to the
directions perpendicular to $\be$; in our case the $(r,\theta)$ directions.
We then use the additional divergence-free condition \eqref{eq:comp_nodiv} to
show that $\kappa$ is harmonic:
\begin{equation}
\left.\begin{array}{r}
\bg=\gradh \kappa \\
\divv \bg = 0 \end{array}\right\}\Rightarrow 
\laph \kappa = 0
\end{equation}
The additional condition \eqref{eq:comp_norm_proj}
then provides a Neumann boundary condition on $\kappa$:
\begin{equation}
\left.\begin{array}{r}
\laph \kappa = 0 \\
\bn\cdot\grad\kappa = 0 \mbox{ on } \pd \domain_h\end{array}\right\}
\Rightarrow \kappa=\kappa_0(e)
\label{eq:equiv_demonstration}
\end{equation}
(see figure \ref{fig:compat}, where $e=z$).
Finally
\begin{equation}
\bg=\grad_h\kappa(e) = 0
\end{equation}
since the $\grad_h$ measures variation in the horizontal directions,
which are perpendicular to the coordinate $e$.
The divergence-free condition \eqref{eq:comp_nodiv} 
is satisfied for $\bu$ since \eqref{eq:def_gu} defines 
$\bg_u$ as a curl. It is satisfied for $\bB$ because \eqref{eq:def_su_sB}
is divergence-free if $\bB$ is, \eg~if $\bB$ is expanded as \eqref{eq:poltor2}.

Condition \eqref{eq:comp_norm_proj}, which is called the compatibility
condition and which ensures the equivalence of both formulations, 
is the projection of the original equations normal to the boundary.  
Its interpretation is quite intuitive: the
compatibility condition preserves information about the original equations
which has been lost by taking the curl.  This procedure is familiar from
simpler contexts: when an equation is differentiated, it must be supplemented
by a constant of integration, which is the evaluation of the original equation
at a point.  Condition \eqref{eq:comp_norm_proj} is sufficient but not unique
-- other boundary conditions ensuring \eqref{eq:equiv_systems} exist.

We can extend the equivalence \eqref{eq:equiv_systems} proved in 
\cite{Marques90} to justify the transformed boundary conditions
often used in practice in the toroidal-poloidal formulation.
To impose the boundary condition $\bu-\bu^{\rm bc}=0$
on a simply connected boundary $\pd\domain$ with normal $\be$,
we substitute for the three vector components the conditions:
\begin{subequations}
\label{eq:equiv_systemsbc}
\begin{alignat}{2}
\be\cdot(\bu-\bu^{\rm bc}) &= 0&&\qquad\text{on $\bdy$}\label{eq:eg0bc}\\
\be\cdot\curl(\bu-\bu^{\rm bc}) &= 0&&\qquad\text{on $\bdy$}\label{eq:ecg0bc}\\
\divh(\bu-\bu^{\rm bc}) &= 0&&\qquad\text{on $\bdy$}\label{eq:comp_nodivbc}\\
\bn\cdot(\bu-\bu^{\rm bc}) &= 0 &&\qquad\text{on $\pd(\bdy)$}\label{eq:comp_norm_projbc}
\intertext{where $\pd(\bdy)$ is the one-dimensional boundary of $\bdy$ and
$\bn$ is perpendicular both to this boundary and to $\be$.
Equation \eqref{eq:comp_nodivbc} can be replaced by}
\divv [(\be\cdot \bu)\be] &= \divv\bu - \divh\bu^{\rm bc} =0
&&\qquad\text{on $\bdy$}
\hspace*{5.5cm}
(\ref{eq:comp_nodivbc}\p)
\nonumber
\end{alignat}
\end{subequations}
where the second equality is valid 
when $\bu$ is divergence-free and the boundary conditions are
homogeneous. 

The transformed boundary conditions \eqref{eq:equiv_systemsbc}
are familiar in the context of a spherical or infinite planar surface,
where the additional condition \eqref{eq:comp_norm_projbc}
is not needed since these surfaces have no boundaries.
For example, in the case of flow between two stationary
infinite planes at $z=\pm 1$, boundary conditions \eqref{eq:equiv_systemsbc} 
take the form:
\begin{subequations}
\label{eq:plane}
\begin{alignat}{2}
w &= 0&&\qquad\text{at $z=\pm 1$}\\
\eta &= 0&&\qquad\text{at $z=\pm 1$}\\
\pdz w &= 0&&\qquad\text{at $z=\pm 1$}
\end{alignat}
\end{subequations}
where $w$ and $\eta$ are the vertical velocity and vorticity.

The derivation of both the gauge and the compatibility conditions depend on
properties of the horizontal Laplace equation; see \eqref{eq:phi_gauge} and
\eqref{eq:equiv_demonstration}.  If the only harmonic function is a constant,
the Neumann condition in \eqref{eq:equiv_demonstration} is superfluous.  This
is the case, for example, on the surface of a sphere. All horizontal
directions are periodic, so that functions which grow monotonically in these
directions are excluded; the compatibility condition \eqref{eq:comp_norm_proj}
can simply be dropped.  However, in domains with more complicated topologies,
such as those bounded by two infinite planes or cylinders considered to be
doubly periodic, additional conditions are necessary in order for
\eqref{eq:curlgrad} to hold. A derivation of these conditions for a general
domain can be found in \cite{Marques90}.

\section{Conditions on the velocity field}
\label{sec:hydro_BC}

We now turn to the conditions to be imposed on the velocity field
in the finite-cylindrical geometry for which $\be=\ez$; see \cite{Marques93}.
Because the next two sections will refer exclusively to the velocity,
we drop the subscript $u$. 
For reference, we write for $\bF$ defined in \eqref{eq:poltor2} 
the identities:
\begin{subequations}
\begin{align}
\bF &= -\ez\times\gradh\psi + \gradh\pdz\phi -\ez\laph\phi\label{eq:bF}\\[.1cm]
\curl\bF &= -\ez\times\gradh\lap\phi +\gradh\pdz\psi -\ez\laph\psi
\label{eq:curlbF}\\[.1cm]
\blap\bF &= -\ez\times\gradh\lap\psi+\gradh\pdz\lap\phi-\ez\laph\lap\phi
\label{eq:lapbF}
\\[.1cm]
\curl\blap\bF &= \ez\times\gradh\lap\lap\phi + \gradh\pdz\lap\psi 
- \ez\laph\lap\psi \label{eq:curlapbF}
\end{align}
\label{eq:facilitate}
\end{subequations}
which will facilitate calculations of vector quantities.
\subsection{Gauge and boundary conditions}
The governing equations are:
\begin{subequations}
\label{eq:potMHD_u_bis}
\begin{align}
(\pdt-\iR\lap)\laph\psi &= \ez\cdot\curl \bsource
\label{eq:potMHD_u1_bis}\\
(\pdt-\iR\lap)\lap\laph\phi &= -\ez\cdot\curl\curl \bsource
\label{eq:potMHD_u2_bis}
\end{align}
\end{subequations}
The system \eqref{eq:potMHD_u_bis} contains five Laplacians
acting in the horizontal directions and three acting in the vertical
directions. Three conditions in each direction are derived from the 
velocity boundary conditions. The two remaining conditions in the
horizontal direction are the gauge and compatibility conditions.

The simplest choice of gauge is:
\begin{subequations}
\begin{alignat}{2}
\phi&=0 &\quad\text{at}\quad & r=1 
\label{eq:gauge_phi_u}
\intertext{along with}
\psi&=0 &\quad\text{at}\quad & r=0
\label{eq:gauge_psi_u}
\end{alignat}
\label{eq:gauge_u}
\end{subequations}%
On the cylinder, boundary conditions are imposed on
$u_r$, $u_\theta$, $u_z$. Referring to \eqref{eq:bF}, we have:
\begin{subequations}
\label{eq:pot_cond_r}
\begin{empheq}[right={\quad\empheqrbrace\quad\text{at}\quad r=1}]{align}
u_r=\ir\pdth\psi+\pdrz\phi&=0\label{eq:pot_cond_r_r}\\
u_\theta=-\pdr\psi+\ir\pdthz\phi&=0\label{eq:pot_cond_r_th}\\
u_z = -\laph\phi &= 0\label{eq:pot_cond_r_z}
\end{empheq}
\end{subequations}
The gauge condition \eqref{eq:gauge_phi_u} 
can be used to simplify \eqref{eq:pot_cond_r_th}:
\begin{equation}
\addtocounter{equation}{-1}
\renewcommand{\theequation}{\thesection.\arabic{equation}b'}
\phi=0 \quad\Rightarrow\quad \pdth\phi
=\pdz\phi=0 \quad\Rightarrow\quad\pdr\psi=0 
\qquad \text{at}\quad r=1
\label{eq:pot_cond_r_th_prime}
\end{equation}

On the (simply-connected) disks, we impose the boundary conditions 
in the form \eqref{eq:equiv_systemsbc}, i.e.
\begin{subequations}
\label{eq:pot_cond_z}
\begin{empheq}[right={\quad\empheqrbrace\quad\text{at}\quad z=\pmh}]{align}
%
0 &= u_z=-\laph \phi\label{eq:eg0bcu}\\
0 &= \ez\cdot\curl\bu = -\laph\psi-\ir\pdr(r^2\angvel_\pm)\label{eq:ecg0bcu}\\ 
0 &= \pdz u_z = -\pdz\laph\phi\label{eq:comp_nodivbcu}
\end{empheq}
\end{subequations}
These are equivalent to those on the individual components but
easier to implement since each of \eqref{eq:pot_cond_z}
involves only one of the potentials.
The remaining condition \eqref{eq:comp_norm_projbc} 
required on the two circles is insured
by \eqref{eq:pot_cond_r_r}.

\subsection{Compatibility condition}
We now turn to the compatibility condition \eqref{eq:comp_norm_proj} for
the hydrodynamic problem in our potential formulation, where
$\be\equiv\ez$, $\bn\equiv\er$, and $\pd\domain_h$ is the $r=1$ boundary:
\begin{equation}
0=\er\cdot\bg=\er\cdot\curl{\mathbf f}
=\er\cdot\curl\left(\left(\pdt-\iR\lap\right)\bu + \bsource\right)
\qquad\text{at}\quad r=1
\label{eq:comp_norm_proj_long}
\end{equation}
Because $\er\cdot\curl$ involves only $\pdth$ and $\pdz$,
derivatives parallel to the $r=1$ boundary, it vanishes 
for all terms in ${\mathbf f}$ which are zero 
or constant at this boundary.
For homogeneous boundary conditions \eqref{eq:pot_cond_r}
on the outer cylinder, this is true for $\pdt\bu$ and for $\bsource$ defined in
\eqref{eq:def_su}
in the absence of a magnetic field, leaving only the Laplacian term.
Referring to \eqref{eq:curlapbF}, we have
\begin{equation}
\er\cdot\curl\blap\bu=\pdrz\lap\psi-\ir\pdth\lap\lap\phi
\end{equation}
Conditions \eqref{eq:gauge_phi_u}, \eqref{eq:pot_cond_r_z}
and (\ref{eq:pot_cond_r_th_prime}) allow the replacement of
$\lap\psi$ and $\lap\phi$ at $r=1$ 
by $\laph\psi$ and $\laph\phi$, which already appear in 
the governing equations \eqref{eq:potMHD_u_bis}, leading to
\begin{equation}
0=\pdrz\laph\psi-\ir\pdth\lap\laph\phi\qquad\text{at}\quad{r=1}
\label{eq:compatibility}
\end{equation}

The complete set of conditions to be imposed on the velocity is then:
\begin{subequations}
\begin{empheq}[right={\quad\empheqrbrace\quad\text{at}\quad r=1}]{align}
\ir\pdth\psi+\pdrz\phi&=0\label{eq:pot_cond_r_r_bis}\\
\pdr\psi&=0 \label{eq:pot_cond_r_th_bis}\\
\laph\phi&=0\label{eq:pot_cond_r_z_bis}\\
\phi&=0 \label{eq:gauge_phi_u_bis}\\
\text{non-axi:}\quad\pdrz\laph\psi-\ir\pdth\lap\laph\phi&=0
\label{eq:compatibility_bis}
\end{empheq}
\label{eq:non_slip_r}
\vspace{-6mm}
\begin{equation}
\hspace{8mm}
\text{axi:}\hspace{36mm}\psi=0\hspace{12.5mm}{\rm at}\quad r=0
\label{eq:gauge_psi_u_bis}
\end{equation}
\end{subequations}
\vspace{1mm} 
\begin{subequations}
\begin{empheq}[right={\quad\empheqrbrace\quad\text{at}\quad z=\pm \frac{h}{2}}]{align}
\hspace{2.6cm}\laph\psi &= -\ir\pdr\left(\angvel_\pm r^2\right)
\label{eq:ecg0bcu_bis}\\ 
\pdz\laph\phi &= 0 \label{eq:comp_nodivbcu_bis}\\
\laph\phi&=0\label{eq:eg0bcu_bis}
\end{empheq}
\label{eq:non_slip_z}
\end{subequations}
These conditions are imposed on $\phi$ and $\psi$ via the influence matrix
method, as will be explained in section \ref{sec:nested_u}.

In equations \eqref{eq:non_slip_r}, we have marked 
conditions \eqref{eq:gauge_phi_u_bis} and 
\eqref{eq:compatibility_bis} as applying 
only to axisymmetric or to non-axisymmetric modes.
This will be explained in the following section.

\subsection{Spatial discretization and symmetry}
\label{sec:axi}

We use the spectral spatial discretization:
\begin{equation}
f(r,\theta,z)=\sum_{m=-\halfM}^{\halfM}f^m(r,z)
e^{im\theta}=\sum_{m=-\halfM}^{\halfM}\,\sum_{k=0}^{K-1}\,
\sum^{2N-1}_
{\genfrac{}{}{0pt}{}{n=|m|} {n+m\text{ even}}}
f_{kn}^m
\mQ_n^m(r)\mT_k\left(\frac{2z}{h}\right)e^{im\theta}
\label{eq:fourdecomp}
\end{equation}
and similarly for $\phi$.  The basis functions in the axial direction $z$ are
the standard Chebyshev polynomials $\mT_k(2z/h)$.  Those in the radial
direction $r$ are the non-standard polynomial basis $\mQ_n^m(r)$ developed by
Matsushima and Marcus \cite{Matsushima}.
Their principal property is that $\mQ_n^m(r) \sim r^m$ as $r\rightarrow 0$,
insuring their regularity at the origin.  The basis functions in the azimuthal
direction $\theta$ are the Fourier modes $e^{im\theta}$.  In
\eqref{eq:fourdecomp}, we do not introduce new notation for Fourier
coefficients, or for coefficients in the 3D tensor-product basis, 
instead distinguishing between physical space values and spectral 
space coefficients by the number and type of superscripts and subscripts.

The decomposition \eqref{eq:fourdecomp} leads to problems and boundary
conditions which are decoupled for each Fourier wavenumber $m$.  In fact,
because of the reflection symmetry in $z$, the problems can be further
reduced.  A vector field is reflection-symmetric in $z$ if $u_r$, $u_\theta$
are even in $z$ and $u_z$ is odd in $z$, \ie~if the potential $\psi$ is even
in $z$ and the potential $\phi$ is odd in $z$, as can be seen from
\eqref{eq:facilitate}.  We denote these functions as having parity
$p=s$.  Quantities related to anti-reflection-symmetric vector fields,
\ie~with $\psi$ odd and $\phi$ even, will be denoted as having parity $p=a$.
The boundary conditions \eqref{eq:non_slip_r} at $r=1$ can be considered
as applying separately to fields of each parity;
note that \eqref{eq:pot_cond_r_r_bis} and \eqref{eq:compatibility_bis}
couple potentials of the same parity.  
The conditions \eqref{eq:non_slip_z} at $z=\pm h/2$ can be reformulated to 
apply only to fields of a single parity; 
for example $\laph\phi(z=\pmh)=0$ can be 
rewritten as $\laph\phi(z=h/2)\pm\laph\phi(z=-h/2)=0$.
Essentiallly, a problem posed over the
entire cylinder can be viewed as $2(M+1)$ problems in the two-dimensional
meridional half-slice $0\leq r \leq 1, 0\leq z \leq h/2$.

Special conditions are applied to the axisymmetric modes.
The gauge freedom \eqref{eq:psi_gauge} for $\psi$
requires the specification of a single value of $\psi$ at each $z$.
In \eqref{eq:gauge_psi_u_bis},
we have chosen to specify this value at the origin:
\begin{equation}
\psi(r=0,\theta,z)
=\sum_m \psi^m(r=0,z)\; e^{im\theta} = \psi^{m=0}(0,z)
\label{eq:phi0}\end{equation}
Condition \eqref{eq:gauge_psi_u_bis} is applied only to the 
axisymmetric mode, since only this mode contributes to the 
sum \eqref{eq:phi0}.

For the axisymmetric modes, two important consequences are derived from the
calculation for an arbitrary function $f^{m=0}(r)$
\begin{equation}
(r\pdr f^0)(r=R) = 
\left. r\pdr f^0\right|_{r=0}^{r=R} 
= \int_0^R \;dr \;\pdr r\pdr f^0 =
\int_0^R r \;dr \;\ir\pdr r\pdr f^0=\int_0^R r \;dr \;\laph f^0
\label{eq:Neu_solvability}\end{equation}
Going from left to right in \eqref{eq:Neu_solvability}, one obtains the
classic solvability condition required by Neumann boundary conditions, since
setting the value of $\pdr f^0(r=R)$ is equivalent to an integral constraint
on $\laph f^0$.  In particular, the Neumann boundary condition
\eqref{eq:pot_cond_r_th_bis} on $\psi$ must be replaced by the integral
constraint like \eqref{eq:Neu_solvability} for the axisymmetric mode.  Going
from right to left in \eqref{eq:Neu_solvability} leads to the conclusion that
the only axisymmetric harmonic function on a disk that includes the origin is
a constant, since $\laph f^0=0$ over $[0,R]$ implies $\pdr f^0(r=R) = 0$ for
each $R$.  This implies that the Neumann boundary condition in
\eqref{eq:equiv_demonstration} is unnecessary to guarantee $\bg=0$ for the
axisymmetric mode, and hence that the compatibility condition
\eqref{eq:compatibility_bis} should not be imposed on the axisymmetric mode.

Corresponding to the removal of the compatibility condition, 
Marques \cite{Marques90}
showed that the system of equations governing the axisymmetric modes
is of lower order. The calculation
\begin{equation}
\pdr^+ f^0 \equiv \frac{1}{r}\pdr r f^0 = 0 \Longrightarrow f^0 = \frac{c}{r} 
\Longrightarrow f^0 = 0
\label{eq:invdrp}\end{equation}
demonstrates the invertibility of $\pdr^+$, or equivalently,
the impossibility on a disk of a non-zero divergence-free axisymmetric
radial vector field which is regular at the origin. Using \eqref{eq:invdrp},
equations \eqref{eq:potMHD_u}-\eqref{eq:potMHD_B} become PDEs of
lower order in $\pdr\psi^0$ and $\pdr\phi^0$:
\begin{subequations}
\begin{align}
\left(\pdt-\iR\lap^+\right)\pdr\psi_u^0 &= \et\cdot\bsource_u^0\\
\left(\pdt-\iR\lap^+\right)\lap^+\pdr \phi_u^0 &= -\et\cdot(\curl\bsource_u^0)\\
\left(\pd_t-\iRm\lap^+\right)\pdr\phi_B^0 &= \et\cdot \bsource_B^0\\
\left(\pd_t-\iRm\lap^+\right)\pdr \psi_B^0 &= \et\cdot\curl \bsource_B^0
\end{align}
\end{subequations}
where $\lap^+\equiv \pdr \pd_r^+ + \pdz^2$.
In the interests of uniformity we continue to solve the
same equations for the axisymmetric as for the non-axisymmetric modes,
altering only the boundary conditions.

\section{Nested Helmholtz and Poisson solvers}
\label{sec:nested_u}

\subsection{Temporal discretization}
\label{sec:discretization_short}

We briefly mention some aspects of our temporal discretization.  A more
extensive description of both the temporal and the spatial discretization is
given in \cite{BoronskiPhD,otherpaper}.
We recall the equations governing the velocity potentials:
\begin{subequations}
\label{eq:potMHD_u_bis2}
\begin{alignat}{3}
(\pdt-\iR\lap)\laph\psi &= \ez\cdot\curl \bsource &&\equiv \source_\psi
\label{eq:potMHD_u1_bis2}\\
(\pdt-\iR\lap)\lap\laph\phi &= -\ez\cdot\curl\curl \bsource &&\equiv \source_\phi
\label{eq:potMHD_u2_bis2}
\end{alignat}
\end{subequations}
Evolution equations such as \eqref{eq:potMHD_u_bis2} are typically
discretized in time via an implicit scheme for the
diffusive terms and an explicit scheme for the nonlinear terms.
For example, with the simplest choice of the backwards and forwards first-order
Euler formulas, the diffusion equation
\begin{subequations}
\label{eq:time_disc}
\begin{align}
\left(\pdt - \iR\lap\right) f &= \source
\label{eq:evolsimp}
\intertext{becomes}
\frac{f(t+\dt) - f(t)}{\dt}-\iR\lap f(t+\dt)&= \source(t)\nonumber\\
\left(I-\frac{\dt}{Re}\lap\right) f(t+\dt) 
&= f(t)+\dt \: \source(t)
\label{eq:discsimp}\end{align}
\end{subequations}
Thus implicit-diffusive/explicit-nonlinear temporal discretization transforms
the pa\-ra\-bo\-lic equation \eqref{eq:evolsimp} into the Helmholtz problem
\eqref{eq:discsimp} for $f(t+\dt)$.  
Similarly, the temporally discretized versions of the more complicated
equations \eqref{eq:potMHD_u_bis2} give $\psi(t+\dt)$ and 
$\phi(t+\dt)$ as solutions to a
sequence of nested Helmholtz and Poisson problems.
We will not distinguish
between the continuous-time parabolic operators of type \eqref{eq:evolsimp}
and the discretized Helmholtz operators of type \eqref{eq:discsimp} and refer
to both as Helmholtz problems.

\subsection{Substitution of Dirichlet boundary conditions}
\label{sec:influence}

Equations \eqref{eq:non_slip_r}-\eqref{eq:non_slip_z} give the set of boundary
conditions which is to be imposed on \eqref{eq:potMHD_u_bis2}.  The major
difficulty of the poloidal-toroidal formulation is that set
\eqref{eq:non_slip_r}-\eqref{eq:non_slip_z}, while appropriate for the entire
problem, does not provide separate boundary conditions appropriate to each
individual Helmholtz and Poisson problem.  Some of the conditions involve both
$\psi$ and $\phi$.  Even conditions involving only one potential can be
problematic because the order of the equations and of the boundary conditions
do not match. The prototypical example of this occurs in the 2D
streamfunction-vorticity formulation.  At each timestep, one would like to
solve successively the Helmholtz problem for the vorticity, and the Poisson
problem for the streamfunction.  However, no boundary conditions are available
for the vorticity, while the streamfunction must satisfy both Dirichlet and
Neumann boundary conditions.

The influence matrix technique \cite{Tuckerman89} calls for replacing the
problematic boundary conditions by conditions which are easier to implement
numerically, in this case Dirichlet boundary conditions on a set of
intermediate fields.  The values used in these boundary conditions are
determined in such a way that the exact boundary conditions are satisfied.  We
show below the sequence of problems with their associated boundary conditions:
\begin{subequations}
\begin{alignat}{3}
&\qquad \left(\pdt-\iR\lap\right)f_\psi = \source_\psi\\[0.1cm]
&\qquad f_\psi = -\ir\pdr(r^2\angvel_\pm)&&\text{at}\quad z=\pmh \label{eq:bcfpsi}\\[-.1cm]
\hspace{-2cm}{\rm axi:} &\qquad\int_0^r r\;dr f_\psi = 0&&\text{at}\quad r=1 \\[-.1cm]
\hspace{-2cm}{\rm nonaxi:}&\qquad 
f_\psi = \sigma_f(z)&&\text{at}\quad r=1\nonumber\\[-.1cm]
&\hspace{2cm}\Uparrow\\[-0.3cm]
&\qquad\ \ 
\truecond_f(z)\equiv \left(\pdrz f_\psi-\ir\pdth g_\phi\right)\biggr\vert_{r=1} &&= 0 
\label{eq:bfdef}\nonumber\\[-.1cm]
\nonumber\\
&\qquad\laph \psi = f_\psi \\[.1cm]
\hspace{-2cm}{\rm axi:}&\qquad\psi = 0 &&\text{at}\quad r=0\\[-.1cm]
\hspace{-2cm}{\rm nonaxi:}&\qquad\pdr\psi = 0 &&\text{at}\quad r=1 
\end{alignat}
\begin{alignat}{4}
&\left(\pd_t-\iR\lap\right)g_\phi = \source_\phi \\[.1cm]
& 
g_\phi=\sigma_g^\pm(r) && \text{at}\quad{z=\pmh}
\nonumber\\[-0.3cm]
&\hspace{1.2cm}\Uparrow\\[-0.1cm]
&\hspace{0.97cm}
\truecond^\pm_g(r)\equiv \pdz f_\phi|_{z=\pmh} = 0 
\label{eq:bgpmdef}\nonumber\\[.1cm]
&
g_\phi = \sigma_g(z)&& \text{at}\quad{r=1}\nonumber\\[-0.2cm]
&\hspace{1.2cm}\Uparrow\\[-0.4cm]
&\hspace{0.97cm} 
\truecond_g(z)\equiv \left(\ir\pdth\psi+\pdrz\phi\right)\biggr\vert_{r=1}&=0 
\nonumber\label{eq:bgdef}\\[-.1cm]
\nonumber\\
&\lap f_\phi = g_\phi \\[.1cm]
&f_\phi=0 &&\text{at}\quad{r=1}\\[-.1cm]
&f_\phi=0 &&\text{at}\quad{z=\pmh}\\[-.1cm]
\nonumber\\
&\laph \phi = f_\phi &\qquad\\[.1cm]
&\phi=0 &&\text{at}\quad{r=1}
\end{alignat}
\label{eq:big_tableau}
\end{subequations}

We have introduced intermediate variables $f_\psi$, $g_\phi$ and $f_\phi$, and
required them to obey Dirichlet boundary conditions with unknown values
$\sigma_f(z)$, $\sigma_g(z)$ and $\sigma^\pm_g(r)$, or 0.  We have also
introduced the notation $\truecond_f(z)$, $\truecond_g(z)$, and
$\truecond^\pm_g(r)$ for quantities which should be zero if the actual
boundary conditions were satisfied.  The boundary conditions in
\eqref{eq:big_tableau} are identical to
\eqref{eq:non_slip_r}--\eqref{eq:non_slip_z}, restated where possible in terms
of $f_\psi$, $g_\phi$ and $f_\phi$.  The influence matrix establishes the
correspondence between $\{\sigma_f,\sigma_g,\sigma^\pm_g\}$ and
$\{\truecond_f,\truecond_g,\truecond^\pm_g\}$. No significance should be
attached to the choice of equation in \eqref{eq:big_tableau} at which each
$\truecond$ is defined, \ie~the elliptic problem with which each of the
original boundary conditions has been associated.  The correspondence serves
merely to establish that the number of unknown Dirichlet values $\sigma$ is
the same as the number of boundary conditions $\truecond=0$.

In order to simplify the notation, we have suppressed the indices labelling
the azimuthal Fourier wavenumber $m$ and axial parity $p\in\{s,a\}$.  Each
equation in \eqref{eq:big_tableau} should in fact be interpreted as applying
separately to modes with different $(m,p)$ values.  Wherever it occurs,
$\pd_\theta$ should be interpreted as multiplication by $im$, while in
equation \eqref{eq:bcfpsi}, the right-hand-side is axisymmetric and hence
should be interpeted as zero for $m\neq 0$.  Note that the boundary conditions
for the axisymmetric and non-axisymmetric modes differ slightly, as explained
in section \ref{sec:axi}.

\subsection{Influence matrix method}\label{sec:infmat}

System \eqref{eq:big_tableau} is solved by generalizing the standard
decomposition of a linear boundary value problem into particular and
homogeneous problems, in which the boundary conditions or right-hand-side are
set to zero, respectively. Here, the nature of the boundary conditions and the
intermediate solutions to which they are applied are also changed. 
Historically, the name capacitance matrix has also been used to
denote what we call the influence matrix; this has guided our 
choice of notation $\infmat$.
The steps for carrying out the influence matrix technique are as follows.

{\bf Preprocessing step (homogeneous solutions):}
\begin{itemize}
\item
We calculate solutions to the homogeneous problem ($\bsource=0$) with a complete
set of Dirichlet boundary conditions corresponding to the spectral
discretization \eqref{eq:fourdecomp}.  Specifically, for each Fourier mode
$m\in\{0,\ldots,\frac{M}{2}\}$ and axial parity $p\in\{s,a\}$, the boundary
values $\{\sigma_f(z),\ \sigma_g(z),\ \sigma^\pm_g(r)\}$ are set successively to:
\begin{subequations}
\begin{alignat}{4}
\Big\{\sigma_f(z)&=\mathcal{T}_k\left(\textstyle\frac{2z}{h}\right), &\quad\sigma_g(z)&=0, 
&\quad\sigma^+_g(r)&=0, &\quad\sigma^-_g(r)&=0\Big\}\label{eq:bcf}\\
\Big\{\sigma_f(z)&=0, &\quad\sigma_g(z)&=\mathcal{T}_k\left(\textstyle\frac{2z}{h}\right), 
&\quad\sigma^+_g(r)&=0, &\quad\sigma^-_g(r)&=0\Big\}\label{eq:bcg}\\
\Big\{\sigma_f(z)&=0, &\quad\sigma_g(z)&=0, 
&\quad\sigma^+_g(r)&=\mathcal{Q}_n^m(r), 
&\quad\sigma^-_g(r)&=\mathcal{Q}_n^m(r)\Big\}\label{eq:bce}\\
\Big\{\sigma_f(z)&=0, &\quad\sigma_g(z)&=0, 
&\quad\sigma^+_g(r)&=\mathcal{Q}_n^m(r),
&\quad\sigma^-_g(r)&=-\mathcal{Q}_n^m(r)\Big\}\label{eq:bco}
\end{alignat}
\end{subequations}

\item
For each homogeneous solution, the values 
of the unsatisfied conditions $\truecond_f$, $\truecond_g$, $\truecond^\pm_g$
are calculated on the boundary. 

\item
These are collected to form the $2(M+1)$ influence matrices $\infmat^{mp}$,
each of size $(K+N)\times(K+N)$. Each set of Dirichlet boundary values
leads to one column of $\infmat^{mp}$.

\item
The influence matrices are inverted to form $(\infmat^{mp})^{-1}$.
Difficulties and techniques related to this inversion are
discussed in appendix \ref{sec:regular_IM}.

\end{itemize}

{\bf Each timestep (particular and final solutions):}

\begin{itemize}
\item
We calculate the particular solution, \ie~the solution to the inhomogeneous
problem ($\bsource\neq 0$) with homogeneous Dirichlet boundary conditions
$\sigma_f=\sigma_g=\sigma^\pm_g=0$.

\item
We calculate the values of the unsatisfied conditions 
$\{\truecond_f,\truecond_g,\truecond^\pm_g\}$ on the boundary.
These are separated according to Fourier mode $m$ and axial parity $p$.

\item
Each set $(m,p)$ of $\truecond$ values is multiplied by the corresponding
matrix $(\infmat^{mp})^{-1}$ to obtain appropriate values of
$\{\sigma_f,\sigma_g,\sigma^\pm_g\}$.

\item
The inhomogeneous problem is then solved again with 
corrected inhomogeneous Dirichlet boundary values.
\end{itemize}

Axial symmetry is taken into account by examining \eqref{eq:big_tableau}.
For solutions which are reflection-symmetric in $z$ ($p=s$), 
$\sigma_f$ is even in $z$, \ie~$k$ takes only even values in 
\eqref{eq:bcf}, while $\sigma_g$ is odd in $z$, so $k$ takes
only odd values in \eqref{eq:bcg}. Additionally, only \eqref{eq:bco}
is used. The corresponding $\truecond_f$ is odd, and $\truecond_g$, $\truecond_g^{pm}$ are even in $z$.
The opposite holds for fields which
are anti-reflection symmetric in $z$ ($p=a$):
$k$ takes only odd values in \eqref{eq:bcf} and even values in 
\eqref{eq:bcg} and only \eqref{eq:bce} is used.

This decomposition can be expressed mathematically as a version 
of the Sherman-Morrison-Woodbury formula \cite{Tuckerman89}.
Here, a large problem coupling $f_\psi,\psi,g_\phi,f_\phi,\phi$ 
is decoupled by a transformation (the change in boundary conditions)
of low rank ($K+N$). The solution to the coupled problem
can be obtained from that of the decoupled problem using an
additional multiplication by a matrix of dimension $K+N$.

\section{Towards an MHD solver}\label{sec:MagneticSolver}

We now address the solution of the induction equation in a finite cylinder of
finite conductivity surrounded by a vacuum extending to infinity. For a sphere
or axially infinite cylinder, the boundary between interior and exterior
domains is associated only with the radial coordinate.  The boundary
surrounding a finite cylinder, however, is specified as a relation between $r$
and $z$.  One approach is to define the induction equation in an integral
formulation. The most important advantage is that no boundary conditions must
be specified. Using this formulation \cite{Dobler98}, a stationary
kinematic dynamo problem was solved in a cylindrical geometry.  In
\cite{Iskakov04,Iskakov05}, a finite volume method is used to discretize the
solution in the interior, which is matched to that in the exterior vacuum via
a boundary element method.  An integral equation
formulation was applied to the entire domain in \cite{Xu04},
and \cite{Guermond07} uses finite elements
with a penalty method to apply boundary conditions. 
To our knowledge, however, there exists as yet no
method applicable to the spectral formulation in a finite cylinder.

\subsection{Matching conditions and gauge}

In the remainder of this section, fields or potentials without subscripts or
superscripts will be taken to refer to the interior magnetic field,
while fields or potentials relating to the field in the exterior vacuum
will be designated by a superscript, \eg~$\bB^{vac}$, $\phi^{vac}$.
We recall the equations governing the interior magnetic potentials:
\begin{subequations}
\label{eq:potMHD_B_bis}
\begin{alignat}{3}
(\pdt-Rm^{-1}\lap)\laph\phi &= \be\cdot \bsource &&\equiv \source_\phi
\label{eq:potMHD_B1_bis}\\
(\pdt-Rm^{-1}\lap)\laph \psi &= \be\cdot\curl &&\bsource \equiv \source_\psi
\label{eq:potMHD_B2_bis}
\end{alignat}
\end{subequations}
System \eqref{eq:potMHD_B_bis} requires two boundary conditions at each
bounding surface and supplementary gauge and compatibility
conditions at the horizontal boundary.  The exterior magnetic field in a
vacuum is described by a single harmonic potential which requires one boundary
condition at each bounding surface; see \ref{sec:exterior}.  Thus a total of
three matching conditions must be applied at each bounding surface:
\begin{subequations}
\label{eq:pot_cond_B}
\begin{empheq}[right=\ \empheqrbrace\quad\text{on}\quad\pd\Omega]{align}
0 &= (B_r-B_r^{vac})=\ir\pdth\psi+\pdr(\pd_z\phi-\phi^{vac})\label{eq:pot_cond_B_r}\\
0 &= (B_\theta-B_\theta^{vac})=-\pdr\psi+\ir\pdth(\pdz\phi-\phi^{vac})
\label{eq:pot_cond_B_th}\\
0 &= (B_z-B_z^{vac}) = -\lap\phi +\pdz(\pdz\phi-\phi^{vac})\label{eq:pot_cond_B_z}
\end{empheq}
\end{subequations}

On the bounding cylinder $r=1$, equations \eqref{eq:pot_cond_B_th} and
\eqref{eq:pot_cond_B_z} (but not \eqref{eq:pot_cond_B_r}) can be simplified by
choosing the gauge:
\begin{equation}
0=(\pdz\phi-\phi^{vac})\qquad\text{at}\quad r=1
\label{eq:gauge_mag}
\end{equation}%
leading to:
\begin{subequations}
\begin{empheq}[right=\quad\empheqrbrace\ \quad\text{at}\quad {r=1}]{align}
0 &=\ir\pdth\psi+\pdr(\pd_z\phi-\phi^{vac})
\label{eq:pot_cond_Br_r}\\
0 &= \pdr\psi\label{eq:pot_cond_Br_th}\\
0 &= \lap\phi\label{eq:pot_cond_Br_z}
\end{empheq}\label{eq:pot_cond_Br}
\end{subequations}

On the disks $z=\pm h/2$, we use the boundary condition \eqref{eq:ecg0bc},
which becomes:
\begin{subequations}
\label{eq:magbcs}
\begin{equation}
0 = \ez\cdot\curl(\bB-\bB^{vac})
=-\laph\psi \qquad\text{at}\quad z=\pmh
\label{eq:nonaxi_contB_z}
\end{equation}
since the exterior magnetic field is curl-free.
Thus $\psi$ is harmonic on the disks, with 
homogeneous Neumann boundary condition \eqref{eq:pot_cond_Br_th} 
and is therefore constant on each disk.
The matching conditions \eqref{eq:pot_cond_B_r} and 
\eqref{eq:pot_cond_B_th} can be applied at the disks to show that
$\pdz\phi-\phi^{vac}$ is constant on each disk, while the
gauge condition \eqref{eq:gauge_mag} shows that the constant is zero:
\begin{equation}
0=\pdz\phi-\phi^{vac} \qquad\text{at}\quad {z=\pmh}
\label{eq:gauge_mag_z}\end{equation}
The matching conditions at $z=\pmh$ are completed by applying
\eqref{eq:pot_cond_B_z} at the disks:
\begin{equation}
0 = -\lap\phi +\pdz(\pdz\phi-\phi^{vac})
 = -\laph\phi -\pdz\phi^{vac}\qquad\text{at}\quad{z=\pmh}
\label{eq:pot_cond_Bz_z}\end{equation}
\end{subequations}
The final set of gauge and matching conditions to be imposed is
\eqref{eq:gauge_mag}, \eqref{eq:pot_cond_Br} and \eqref{eq:magbcs}.

\subsection{Exterior magnetic field}
\label{sec:exterior}

Since the exterior magnetic field has both zero curl and zero divergence,
and the exterior domain is simply connected, \eqref{eq:curlgrad} states
that the exterior magnetic field can be represented as:
\begin{align}
\bB^{vac}&=\grad\phi^{vac}
\label{eq:Bvac}\end{align}
where $\phi^{vac}$ satisfies:
\begin{subequations}
\begin{align}
\lap\phi^{vac} &= 0  \qquad\text{outside the cylinder}\label{eq:phiharm}\\
\grad\phi^{vac} &= 0 \qquad\text{at infinity}\label{eq:gradphidecay}
\end{align}
\label{eq:Bext}
\end{subequations}
Equation \eqref{eq:gradphidecay} supplies the boundary condition on
$\phi^{vac}$ at infinity, while conditions on the cylindrical
boundary are provided by coupling with the interior potentials via the gauge
and matching conditions \eqref{eq:gauge_mag}, \eqref{eq:pot_cond_Br_r},
\eqref{eq:gauge_mag_z} and \eqref{eq:pot_cond_Bz_z}.  

Although \eqref{eq:Bext} is posed in the infinite domain outside the cylinder,
we can avoid discretizing the infinite domain and solving numerically by using
known analytic solutions to the Laplace equation.  Our approach is to
formulate a complete set of analytic solutions to \eqref{eq:Bext}, each of
whose derivatives can be calculated.  The exterior solution $\phi^{vac}$ can
be expanded in this set, with coefficients related to values on the
cylindrical boundary.  The normal derivatives at the boundary can then be
evaluated in terms of these coefficients.  This defines a correspondence
between a set of boundary values $\{\phi^{vac}|_{\bdy}\}$ and a set of normal
derivatives $\{\bn\cdot\grad\phi^{vac}|_{\bdy}\}$.  This correspondence, or
influence matrix, constitutes a basis for the Dirichlet-to-Neumann mapping for
the domain outside the finite cylinder.  The normal derivatives appearing in
the matching conditions \eqref{eq:pot_cond_Br_r} and \eqref{eq:pot_cond_Bz_z}
can then be replaced by functions $\mF_r\equiv \pdr\phi^{vac}|_{r=1}$ and
$\mF^\pm_z\equiv \pdz\phi^{vac}|_{z=\pmh}$ of the boundary values
$\{\phi^{vac}|_{\bdy}\}$.  Equations \eqref{eq:gauge_mag} and
\eqref{eq:gauge_mag_z} in turn relate the boundary values of the exterior and
interior potentials via $\phi^{vac}|_{\bdy}=\pdz\phi|_{\bdy}$.  The exterior
magnetic field no longer appears and the interior problem is closed.
Essentially, we seek to replace the matching conditions
\eqref{eq:gauge_mag},\eqref{eq:gauge_mag_z},\eqref{eq:pot_cond_Br_r},
\eqref{eq:pot_cond_Bz_z} by:
\begin{subequations}
\label{eq:DtN}
\begin{align}
0 &=\ir\pdth\psi+\pdrz\phi
-\mF_r(\{\pdz\phi|_{\bdy}\}) && \hspace{-2cm}\text{at}\quad r=1 \\
0 &= -\laph\phi-\mF^\pm_z(\{\pdz\phi|_{\bdy}\}) &&\hspace{-2cm}\text{at}\quad z=\pmh
\end{align}
\end{subequations}
The task is now to obtain a well conditioned matrix representation of the
Dirichlet-to-Neumann mappings $\mF_r,\mF^\pm_z$.

We have considered two sets of solutions to \eqref{eq:Bext}.
The first set is constructed from the classic spherical harmonics.
That is, we expand $\phi^{vac}$ as:
\begin{EQ}[rcl]
\phi^{vac}&=&\phi^{vac}_\infty+\sum_m \sum_{l\geq |m|} 
\phi^{vac}_{lm}\rho^{-(l+1)} P_{lm}(\cos\xi)e^{im\theta}\\
\rho&=&\sqrt{r^2+z^2}, \qquad\xi=\tan^{-1}\left(\frac{r}{z}\right)
\label{eq:phi_as_spher_harm}\\
\end{EQ}
where $P_{lm}$ are the associated Legendre polynomials.
According to \eqref{eq:Bvac}, $\phi^{vac}$ is defined only up 
to a constant, which we may choose such as to set $\phi^{vac}_\infty=0$.
(In two dimensions, \ie~for a function whose gradient decays as
the cylindrical radius $r$, rather than the spherical radius $\rho$,
tends to infinity, logarithmic functions would have to be included 
in the expansion because it cannot be assumed that $\phi^{vac}$ tends
to a constant at infinity.)
For a given longitudinal Fourier mode $m$, the solution
\eqref{eq:phi_as_spher_harm} has degrees of freedom associated 
with index $l$, associated with the latitude.
In the standard spectral-physical space duality,
the set of coefficients $\phi^{vac}_{lm}$ corresponds to the
set of values $\phi^{vac}_m(r_i,z_i)$ for $(r_i,z_i)$ 
on the boundary and can be determined from them by solving:
\begin{EQ}[rcl]
\phi^{vac}_m(r_i,z_i)&=&\sum_{l=|m|}^{|m|+L-1}
\phi^{vac}_{lm}\rho_i^{-(l+1)} P_{lm}(\cos\xi_i)\\
\rho_i&=&\sqrt{r_i^2+z_i^2}, \qquad\xi_i=\tan^{-1}\left(\frac{r_i}{z_i}\right)
\label{eq:phi_as_spher_harm_disc}
\end{EQ}
where $L$ is the number of points on the boundary.
Expansion \eqref{eq:phi_as_spher_harm} readily yields the normal
derivatives $\pdr\phi^{vac}_m(r_i,z_i), \pdz\phi^{vac}_m(r_i,z_i)$ 
at the boundary in
terms of the coefficients $\phi^{vac}_{lm}$.
However, this approach is not feasible, in part because the transform
\eqref{eq:phi_as_spher_harm_disc} is extremely poorly conditioned,
like all other transforms involving monomials.
As $l$ increases, the functions $\rho_i^{-(l+1)}$ become
spiked at the largest values of $\rho_i$ and zero elsewhere,
a difficulty which does not arise on a spherical surface
where $\rho$ is constant. 

We have also considered a second set of solutions to \eqref{eq:Bext},
constructed from the equally classic free-space or fundamental 
Green's functions:
\begin{equation}
\phi^{vac}(\bx)=\int_{\bdy}\;\frac{d\bx\p}{4\pi\;|\bx-\bx\p|}\;\sigma(\bx\p)
\label{eq:Green}
\end{equation}
where $\sigma(\bx\p)$ is a distribution on the cylindrical surface $\pd\domain$ 
which is calculated in such a way as to yield a particular set of 
boundary values for $\phi^{vac}(\bx)$ and $G(\bx,\bx\p)=4\pi/|\bx-\bx\p|$
satisfies
\begin{equation}
\lap_\bx\p G(\bx,\bx\p) = \delta(\bx-\bx\p) 
\end{equation}
The derivatives of expansion \eqref{eq:Green} are obtained by
differentiating the Green's functions:
\begin{equation}
\grad_{\bx}\phi^{vac}(\bx)
=\int_{\bdy} d\bx\p\;\frac{\bx-\bx\p}{4\pi\;|\bx-\bx\p|^3}\;\sigma(\bx\p)
\end{equation}
This approach is discussed and shown to perform quite well for a
two-dimensional test problem in \cite{BoronskiJCP}.  The exterior fields
generated approximate the solution uniformly near the boundary and converge
exponentially with resolution; the influence matrix representing the
Dirichlet-to-Neumann mapping is well conditioned.

\subsection{Magnetic compatibility condition}

The magnetic compatibility condition, obtained by substituting 
\eqref{eq:def_gB} into \eqref{eq:comp_norm_proj}, is:
\begin{equation}
0 = \er\cdot\bg_B
=\er\cdot\left(\left(\pdt-Rm^{-1}\lap\right)\bB + \bsource_B\right)\qquad\text{at}\quad r=1
\end{equation}
Contrary to the velocity, the magnetic field does not vanish on the
boundary so $\pdt\bB\neq 0$. The nonlinear term 
$\bsource_B=\curl(\bu\times\bB)$ does vanish at the boundary since 
$\bu|_{r=1}=0$ and its radial curl $\er\cdot\bsource_B$ 
contains no normal derivatives of $\bu$.
The remaining terms are evaluated using \eqref{eq:bF} and \eqref{eq:lapbF}:
\begin{equation}
0=\er\cdot (\pdt-\iRm\blap)\bB = \ir\pdth(\pdt-\iRm\lap)\psi + 
\pdrz(\pdt-\iRm\lap)\phi \quad\text{at}\quad r=1
\label{eq:mag_compat}
\end{equation}
The magnetic compatibility equation must use the same time discretization
as the evolution equations, here backwards Euler.
Although the the time derivative in \eqref{eq:mag_compat} may seem difficult to
include in an implementation of the influence matrix, the
boundary operator in \eqref{eq:mag_compat} can be decomposed into two parts,
one which acts on the homogeneous solution at time $t+\dt$ (and contributes to the influence matrix) and the other 
which acts on the particular solution at time $t+\dt$ and the actual 
solution at time $t$.  Details are given in
\cite{BoronskiPhD}. 

The velocity compatibility condition \eqref{eq:compatibility_bis}
must also be modifed to include a contribution from the magnetic field:
\begin{equation}
0=\pdrz\laph\psi_u-\ir\pdth\lap\laph\phi_u
-\er\cdot\curl(\bB\cdot\bgrad)\bB \qquad\text{at}\quad r=1
\label{eq:vel_comp_with_mag}
\end{equation}

\subsection{Nested elliptic problems and influence matrix}

The set of nested Helmholtz and Poisson problems \eqref{eq:potMHD_B_bis},
together with the conditions 
\eqref{eq:pot_cond_Br_th}-\eqref{eq:pot_cond_Br_z},
\eqref{eq:nonaxi_contB_z}, \eqref{eq:DtN} and \eqref{eq:mag_compat}
can be solved using the influence
matrix technique, as was done for the velocity in section 
\ref{sec:influence} and in equation \eqref{eq:big_tableau}.
\begin{subequations}
\begin{alignat}{3}
&\quad \left(\pdt-\iR\lap\right)f = \source_\psi\label{eq:fpsi_mag}\\
&\quad f = 0 &&\text{at}\quad z=\pmh\\[-.3cm]
{\rm axi:} &\quad\int_0^r r\;dr f = 0 &&\text{at}\quad r=1\\[-.1cm]
{\rm nonaxi:} &\quad f = \sigma_f(z)&&\text{at}\quad r=1\nonumber\\[-0.2cm]
&\hspace{1.4cm}\Uparrow\\[-0.4cm]
&\hspace{1.2cm} \label{eq:compat_B_with_dt}
\truecond_f(z)\equiv \left[\ir\pdth(\pdt-\iRm \lap)\psi\right. 
&&+
\left.\pdrz(\pdt-\iRm \lap)\phi\right]\biggr\vert_{r=1} = 0\nonumber
\\[-.1cm]
\nonumber\\
&\quad\laph \psi = f \\
{\rm axi:}&\quad\psi = 0 &&\text{at}\quad r=0\\[-.1cm]
{\rm nonaxi:}&\quad\pdr\psi = 0 &&\text{at}\quad r=1\\[-.1cm]
\nonumber
\end{alignat}
\begin{alignat}{3}
&\hspace{-0.26cm}\left(\pdt-\iR\lap\right)g = \source_\phi 
\label{eq:gphi_mag}\\[.1cm]
&\hspace{-0.26cm}g=\sigma_g(z)&&\hspace{0.1cm}\text{at}\quad r=1
\nonumber\\[-0.2cm]
&\hspace{0.7cm}\Uparrow\\[-0.2cm]
&\hspace{0.5cm}
\truecond_g(z)\equiv \lap\phi|_{r=1} = 0 \nonumber
\\
& 
\hspace{-0.26cm}g = \sigma^\pm_g(r)&&\hspace{0.1cm}\text{at}\quad z=\pmh
\nonumber\\[-0.2cm]
&\hspace{0.7cm}\Uparrow\\[-0.1cm]
&\hspace{0.5cm} 
\truecond^\pm_g(r)\equiv -g|_{z=\pmh}
&&\hspace{-1.8cm}-\mF^\pm_z(\{\pdz\phi|_{\bdy}\})=0\nonumber
\\[-.1cm]
\nonumber\\
&\hspace{-0.26cm}\laph \phi = g \\[.1cm]
& 
\hspace{-0.26cm}\phi=\sigma_\phi(z)&&\hspace{0.1cm}\text{at}\quad r=1
\nonumber\\[-0.2cm]
&\hspace{0.7cm}\Uparrow\\[-0.4cm]
&\hspace{0.5cm}
\truecond_\phi(z)\equiv \left[\ir\pdth\psi+\pdrz\phi\right]\Bigg\vert_{r=1}
&&\hspace{-0.5cm}-\mF_r(\{\pdz\phi|_{\bdy}\})=0\nonumber
\end{alignat}
\label{eq:big_tableau_mag}
\end{subequations}\\
The major differences between the velocity and magnetic cases
are the reduction from five to four in the number
of elliptic problems, the time derivative in the magnetic
compatibility equation, and the presence of 
the Dirichlet-to-Neumann mappings $\mF_r$ and $\mF^\pm_z$ which have 
replaced the normal derivatives of the exterior solution $\phi^{vac}$.
We recall the meaning of $\mF_r(\{\pdz\phi|_{\bdy}\})$ and 
$\mF^\pm_z(\{\pdz\phi|_{\bdy}\})$:
the set $\{\pdz\phi|_{\bdy}\}$ provides Dirichlet boundary values for
the exterior Laplace problem and $\mF_r$ and $\mF^\pm_z$
are the normal derivatives of the exterior solution at the boundaries.

The Dirichlet-to-Neumann mappings, as well as each solution and problem in
\eqref{eq:big_tableau_mag}, decouple according to azimuthal Fourier mode $m$
and axial parity $p\in\{s,a\}$.  As before, we construct the influence matrix
by solving homogeneous versions of \eqref{eq:big_tableau_mag}, with $\source_\psi$
and $\source_\phi$ set to zero in \eqref{eq:fpsi_mag} and \eqref{eq:gphi_mag} and a
complete set of Dirichlet boundary values $\sigma_f$, $\sigma_g$,
$\sigma^\pm_g$, $\sigma_\phi$.  Evaluating $\truecond_f$, $\truecond_g$, $\truecond^\pm_g$, $\truecond_\phi$
yields the influence matrix.

\section{Conclusion}

We have presented a poloidal-toroidal formulation for solving the
time-dependent three-dimensional magnetohydrodynamic equations in a finite
cylinder.  While preserving the original mathematical formulation described in
\cite{Marques90} and later tested on a linear stability Rayleigh-B\'enard
convection problem in \cite{Marques93}, we incorporated the influence matrix
technique \cite{Tuckerman89} for decoupling the boundary and compatibility
conditions emerging from the potential formulation.  We have also described an
extension of this algorithm to the induction equation governing the evolution
of the magnetic field.

The most important advantage of using the toroidal-poloidal decomposition is
that the divergence-free character of the velocity and magnetic fields is
imposed exactly, by construction.  For the induction equation, the potential
formulation makes it possible to solve for the magnetic field without
introducing an artificial numerical magnetic analogue to the hydrodynamic
pressure which has no physical meaning. In addition, using scalar functions
instead of components of vector fields simplifies and homogenizes the usage of
differential operators.  The influence matrix technique allows the
poloidal-toroidal formulation to be sufficiently economical to be used for
time-integration.

We have im\-ple\-men\-ted and validated this method for the hydrodynamic von
K\'arm\'an problem of flow in a cylinder driven by counter-rotating disks,
using a spectral discretization which is regular on the cylindrical axis.
These results are presented in a companion article \cite{otherpaper}.  In
extending this algorithm to the full magnetohydrodynamic problem, no
difficulty is posed by the induction equation, whose structure is simpler
than that of the Navier-Stokes equation.  Instead, the main difficulty is that
the magnetic field is not specified at the domain boundary but must instead
satisfy matching conditions between the interior domain (here a finite
cylinder) and the exterior domain (here an infinite vacuum).  We have
developed a formalism involving the Dirichlet-to-Neumann mapping for
eliminating the exterior magnetic field, which has been implemented and
validated for a two-dimensional test problem \cite{BoronskiJCP}.  Future
research will focus on implementing the poloidal-toroidal formulation for the
full magnetohydrodynamic problem.

\newcounter{oldfigure}
\setcounter{oldfigure}{\value{figure}}

\appendix

\addtocounter{oldfigure}{1}
\renewcommand{\thefigure}{\arabic{oldfigure}}

\section{Regularizing the influence matrix}
\label{sec:regular_IM}

In section \ref{sec:infmat}, we have assumed that the influence matrices
$\infmat^{pm}$ are invertible, which is actually not the case.  The influence
matrices are non-invertible for several reasons. The first issue is
geometric. The finite cylinder has corners at which conditions at $z=\pmh$ and
$r=1$ must both be satisfied. When formulated in spectral rather than physical
space, the redundant conditions correspond to a linear combination of rows and
cannot be easily identified.  A second factor is the discretization of the
Poisson and Helmholtz solvers, in particular the replacement of the
highest-wavenumber equations by the boundary conditions mandated by the $\tau$
method.  A third cause is the decrease in polynomial order due to
differentiation by boundary operators. For numerical reasons,
the eigenvalues corresponding to these directions may be nearly zero,
rather than exactly so.

One remedy \cite{Tuckerman89} consists of thresholding: diagonalizing
$\infmat^{pm}$ and replacing the eigenvalues whose absolute values are below a
certain experimentally-determined threshold $\epsilon_\mu$ by an arbitrary
value, say 1, leading to an invertible matrix.  The justification of this
manipulation of the spectrum is that the eigenvectors corresponding to the
zero eigenvalues play no role in satisfaction of the boundary conditions.
This is true if the linear system of equations defined by the influence matrix
and the right-hand-side is underdetermined, \ie~if the right-hand-side belongs
to the image space of the influence matrix.  Because the particular solutions
are determined using the same nested solver used for constructing the
homogeneous solutions, this is in fact the case.

However, a major problem remains. Even after eigenvalues are eliminated which
would be exactly zero if infinite precision were used, the resulting matrices
still have very small eigenvalues, i.e.~they are still poorly conditioned.
There are various causes for this.  Some boundary value distributions are
almost linearly dependent.  More importantly, because some boundary conditions
are of higher differential order then others, the magnitudes of different
portions of the influence matrices are very different.  We shall call these
eigenvalues {\it small}, in contrast to those which would be zero in infinite
precision, which we shall call simply zero eigenvalues.  The small eigenvalues
depend on the spatial resolution and on the product $Re/\dt$, a parameter
which appears in the Helmholtz problem.  As the resolution or $Re/\dt$ are
increased, an increasing number of small eigenvalues appear, whereas the
number of zero eigenvalues depends only on the geometry and on the kind of
boundary conditions.  

The condition number, approximately the ratio between the largest and smallest
eigenvalues values of a matrix, is an upper bound on the number of
lost meaningful digits in the numerical solution to a linear equation
involving this matrix.  If the threshold $\epsilon_\mu$ is chosen such as to
lower the condition number to an acceptable value ($O(10^8)-O(10^{10}$), then
small eigenvalues are eliminated, in addition to zero eigenvalues, leading
to errors in satisfaction of the boundary conditions that are much
higher than machine precision.
The challenge is, first, to distinguish between the zero and small
eigenvalues so as to eliminate only the zero eigenvalues and, second, to 
improve the condition number of the adjusted matrix in some way other than
by eliminating the small eigenvalues.

We first modified the thresholding procedure by using the singular value
decomposition (SVD) rather than diagonalization. The advantages of the SVD is,
first, that it always exists and, second, that the matrix of singular vectors
is better conditioned than the eigenvector matrix, because the left and right
singular vectors are orthogonal, in contrast to eigenvectors, which may be
close to linearly dependent.  As an example, we consider the influence matrix
for spatial resolution $(N=96)\times(K=192)$, Reynolds number $Re=10^{4}$ and
time step $\dt=10^{-2}$. The magnitudes of the singular values $\gamma_i$ of
the influence matrix block $\infmat^{1,s}$ with azimuthal Fourier wavenumber
$m=1$ and axial parity $p=s$ are presented on figure \ref{fig:svd_sv}a.
$\infmat^{1,s}$ has one zero singular value (i.e.~a singular value which would
be exactly zero in infinite precision).  Figure \ref{fig:svd_sv}(left) shows
that this value is numerically $10^{-21}$, separated from the next smallest
singular value. In contrast, the zero eigenvalue and smallest remaining
eigenvalue are of the same size and hence cannot be distinguished.  However,
thresholding, whether by replacing this singular value or the corresponding
eigenvalue, still leaves the condition number unacceptably high, on the order
of $10^6/10^{-15}=10^{20}$ for this example.

\begin{psfrags}
\psfrag{mu}[r][][1.]{$\gamma_i$}
\psfrag{i}[r][][1.]{$i$}
\psfrag{\ 1e002}[c][][0.625]{$10^2$}
\psfrag{\ 1e000}[c][][0.625]{$1$}
\psfrag{\ 0.01}[c][][0.625]{$10^{-2}\;$}
\psfrag{\ 0.0001}[c][][0.625]{$10^{-4}$}
\psfrag{\ 1e-006}[c][][0.625]{$10^{-6}$}
\psfrag{\ 1e-008}[c][][0.625]{$10^{-8}$}
\psfrag{\ 1e-010}[c][][0.625]{$10^{-10}$}
\psfrag{\ 1e-012}[c][][0.625]{$10^{-12}$}
\psfrag{\ 1e-014}[c][][0.625]{$10^{-14}$}
\psfrag{\ 1e-016}[c][][0.625]{$10^{-16}$}
\psfrag{\ 1e+006}[c][][0.625]{$\;\;\;10^{6}$}
\psfrag{\ 1000}[c][][0.625]{$10^{3}\ $}
\psfrag{\ 0.001}[c][][0.625]{$10^{-3}$}
\psfrag{\ 1e-006}[c][][0.625]{$10^{-6}$}
\psfrag{\ 1e-009}[c][][0.625]{$10^{-9}$}
\psfrag{\ 1e-012}[c][][0.625]{$10^{-12}$}
\psfrag{\ 1e-015}[c][][0.625]{$10^{-15}$}
\psfrag{\ 1e-018}[c][][0.625]{$10^{-18}$}
\psfrag{\ 1e-021}[c][][0.625]{$10^{-21}$}
\begin{figure}[h]
\centering
	\includegraphics[width=0.45\textwidth]{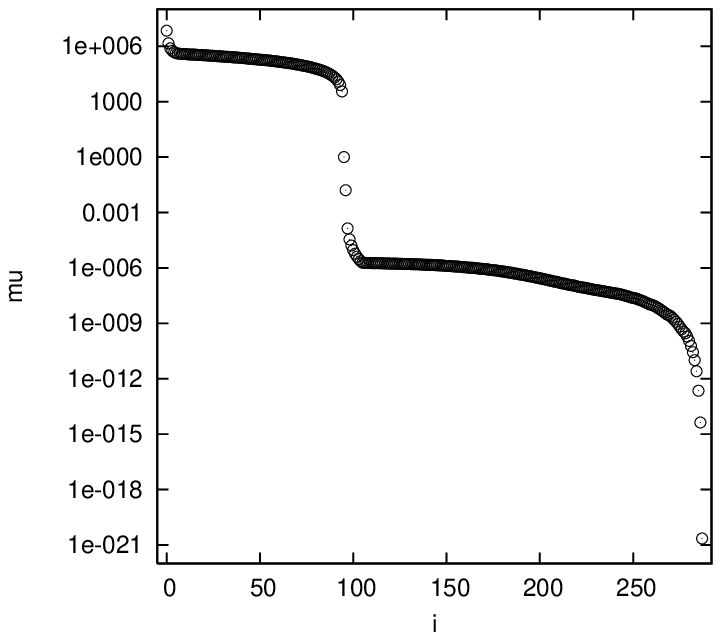}
	\qquad \psfrag{mu}[r][][1.]{$\gamma_i'$}
	\includegraphics[width=0.45\textwidth]{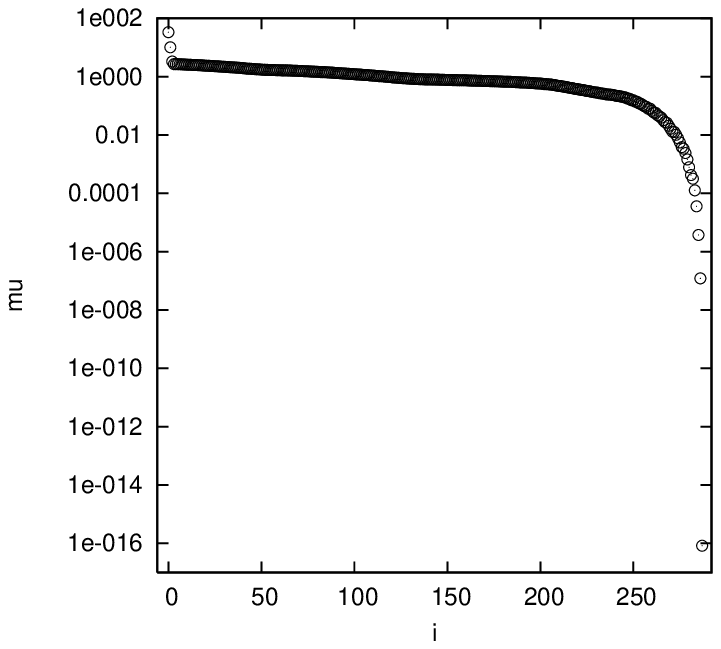}
\caption{Singular values $\{\gamma_i\}$ for case with $(N=96)\times(K=192)$,
$Re/\dt =10^6$, $m=1$, $p=s$. Left: original influence matrix
$\infmat^{1,s}$. Right: scaled influence matrix
$({\infmat^{1,s}})^\prime\equiv
\mx{\boldsymbol{\;\alpha}\;}\infmat^{1,s}\mx{\;\boldsymbol{\beta}\;}$.}
\label{fig:svd_sv}
\end{figure}
\end{psfrags}

The matrix condition number can be decreased more effectively by scaling prior
to thresholding.  If each row is divided by its norm, the condition number of
the matrix is significantly reduced, down to $10^{11}$ for $m>0$ and below
$10^8$ for $m=0$. Additional scaling of columns does not significantly change
the condition number.  The condition number of the $m=0$ matrices are
sufficiently decreased by scaling only their rows, because their poloidal and
toroidal potentials are not coupled.  For $m>0$, the influence matrices can be
further improved by scaling blocks corresponding to different combinations of
types of test functions and boundary conditions.  The influence matrices
$\infmat^{pm}$ for $m>0$ are composed of $9$ submatrices: the three columns
correspond to the three types of imposed simplified Dirichlet boundary
conditions $\{\sigma_g(z),\sigma_f(z),\sigma_g^\pm(r)\}$ while the three rows
correspond to the three quantities
$\{\truecond_g(z),\truecond_f(z),\truecond_g^\pm(r)\}$ reflecting the actual
boundary conditions.  We represent the $|\cdot|_\infty$ norm of the
corresponding submatrix of $\infmat^{pm}$ by $c_{ij}$, where, for
example, $c_{12}$ designates the norm of the $(c_g,\sigma_f)$ submatrix.

For the resolution $(N=96)\times(K=192)$,
the matrix norms $c_{ij}$ of each of the blocks are:
\begin{equation}
	\infmat^{1,s}\quad = \quad
\left[
\begin{BMAT}(e)[1pt,2.8cm,2.8cm]{c.c.c}{c.c.c}
c_{11} & c_{12} & c_{13}\\
c_{21} & c_{22} & c_{23}\\
c_{31} & c_{32} & c_{33}\\
\end{BMAT}
\right]
=\left[
\begin{BMAT}(e)[1pt,2.8cm,2.8cm]{c.c.c}{c.c.c}
10^7 & 1 & 1\\
10^{-2} & 10^{-5} & 10^{-5}\\
0 & 10^{-5} & 10^{-4}\\
\end{BMAT}
\right]
\end{equation}

We now wish to
scale the block-rows and block-columns in such a way as to make the norms 
$c_{ji}'$ of the resulting scaled blocks equal to one another.
\begin{equation}
\mx{\boldsymbol{\alpha}}\infmat^{1,s}\mx{\boldsymbol{\beta}}
\quad = \quad
\left[
\begin{BMAT}(r)[2pt,2.0cm,2.0cm]{c.c.c}{c.c.c}
\alpha_1 & 0 				& 0\\
 0			 & \alpha_2 & 0\\
 0			 & 0				& \alpha_3\\
\end{BMAT}
\right]
\left[
\begin{BMAT}(e)[2pt,2.0cm,2.0cm]{c.c.c}{c.c.c}
c_{11} & c_{12} & c_{13}\\
c_{21} & c_{22} & c_{23}\\
0 		 & c_{32} & c_{33}\\
\end{BMAT}
\right]
\left[
\begin{BMAT}(r)[2pt,2.0cm,2.0cm]{c.c.c}{c.c.c}
\beta_1 & 0 				& 0\\
 0			 & \beta_2 & 0\\
 0			 & 0				& \beta_3\\
\end{BMAT}
\right]
=
\left[
\begin{BMAT}(e)[2pt,2.0cm,2.0cm]{c.c.c}{c.c.c}
c_{11}' & c_{12}' & c_{13}'\\
c_{21}' & c_{22}' & c_{23}'\\
0 		 & c_{32}' & c_{33}'\\
\end{BMAT}
\right]
\label{eq:scaling_matrix}
\end{equation}
In general there exist no $\{\alpha_i\}$ and $\{\beta_i\}$ satisfying
$c_{11}'=c_{12}'=c_{13}'=c_{21}'=c_{22}'=c_{23}'=c_{32}'=c_{33}'=1$. 
We can instead
require
\begin{EQ}[lll]
\alpha_1\beta_1c_{11} &= \alpha_2\beta_2c_{22} &= \alpha_3\beta_3c_{33}=1\\
\alpha_2\beta_1c_{21} &= \alpha_1\beta_2c_{12} &\\
\alpha_3\beta_1c_{31} &= \alpha_1\beta_3c_{13} &
\label{eq:equalizing_cond}
\end{EQ}
The system \eqref{eq:equalizing_cond} has an infinite number of possible
solutions, from which we can select the following:
\begin{EQ}[lll]
\alpha_1=\sqrt{\frac{c_{21}c_{32}c_{33}}{c_{11}c_{12}c_{23}}} & 
\alpha_2=\sqrt{\frac{c_{32}c_{33}}{c_{22}c_{23}}} &
\alpha_3=1 \\
\beta_1=\sqrt{\frac{c_{11}c_{23}}{c_{11}c_{21}c_{32}c_{33}}} & 
\beta_2=\sqrt{\frac{c_{23}}{c_{22}c_{32}c_{33}}} &
\beta_3=\frac{1}{c_{33}} 
\label{eq:scaling_factors}
\end{EQ}
After scaling using \eqref{eq:scaling_matrix}-
\eqref{eq:scaling_factors}, the influence matrix $\infmat^{1,s}$ 
has the following structure:
\begin{equation}
{\infmat^{1,s}}^\prime\quad = \quad
\left[
\begin{BMAT}(e)[1pt,2.8cm,2.8cm]{c.c.c}{c.c.c}
c_{11}' & c_{12}' & c_{13}'\\
c_{21}' & c_{22}' & c_{23}'\\
c_{31}' & c_{32}' & c_{33}'\\
\end{BMAT}
\right]
=
\left[
\begin{BMAT}(e)[1pt,2.8cm,2.8cm]{c.c.c}{c.c.c}
1 & 10^{-2} & 10^{-2}\\
10^{-2} & 1 & 1\\
0 & 1 & 1\\
\end{BMAT}
\right]
\end{equation}
The zero singular value is then easily identified and replaced by 1, leading
to a condition number of $10^{11}$, like that for simple row or column
scaling.  But if the block scaling is followed by row scaling and then
replacement of the zero singular value, then the condition number is further
reduced to the acceptable value of $10^8$.  The singular values of the matrix
after block and row scaling are presented on figure \ref{fig:svd_sv}(right).

We note that operators with high condition numbers are inherent
in the numerical discretization of partial differential equations;
for example, the 1D second derivative operator with homogeneous boundary
conditions using a basis of $K$ Chebyshev polynomials and corresponding grid
has condition number $O(K^4)$. 
The requirements for the solution of the linear systems that occur
in this context are not those of numerical linear algebra:
the right-hand-side is not arbitrary, but results from
time-integration, and not all components are of equal weight.
In our case, we find that after an initial small integration time of 
$T=100\:\dt$, the right-hand-side is always such that the influence matrix,
scaled and regularized to reduce its condition number to $O(10^8)$,
can be inverted to satisfy the constraints in system \eqref{eq:big_tableau} 
to machine accuracy \cite{otherpaper}.

One welcome consequence of scaling is that it separates the zero singular
values which result from non-invertibility of the matrix from the small
singular values which result from poor conditioning of the various components
of the influence matrix).  Without scaling, the number of singular values
below a fixed threshold depends on the spatial resolution and so the zero
eigenvalues cannot be reliably identified and removed. More importantly,
scaling vastly improves the condition number of the influence matrix,
insuring satisfaction of the constraints to machine accuracy.

\end{document}